\documentclass{article}

\date{June 2026}

\title{Paranatural Category Theory}
\author{Jacob Neumann}

\usepackage{agda/agda}
\usepackage{agda/agda-supplement}
\usepackage{sty/packages}
\usepackage{sty/macros}
\usepackage{sty/inclusion}

\begin{document}
\maketitle

\section{Introduction}\label{intro}

\begin{note}
The original ``\texttt{v1}'' version of this article~\cite{neumann2023paranatural} contained several incorrect claims:\footnote{As announced at TYPES 2024~\cite{neumann2024updates}.} see \cref{errata}. The present version is a work-in-progress (to be replaced with a finished version with full discussion), but the incorrect claims have been replaced with correct ones.
\end{note}

\clearpage

\section{Difunctors and Strong Dinaturality}\label{difunctorsStrongDinaturality}

\begin{note}
Throughout, we'll use the term ``difunctor'' to refer to functors of shape $\C\op\times\C\to\SET$. This is more specific than ``profunctor'', where the contravariant and covariant arguments can be different categories, i.e. $\C\op\times\mb D\to\SET$. Many of these results will still hold when the codomain set is some arbitrary category, not $\SET$ specifically.
\end{note}

\mkDefn[\formalLink{StrongDinaturalTransformation}]{StrongDinaturalTransformation}

\begin{figure}
\[
\begin{tikzcd} & \Gamma(I,I) \arrow{rrr}{\alpha_I} \arrow{dr}{\Gamma(I,j)} &&& \Delta(I,I) \arrow{dr}{\Delta(I,j)} \\
\UNIT \arrow{ur}{x} \arrow[swap]{dr}{y} && \Gamma(I,J) &&& \Delta(I,J) \\
 & \Gamma(J,J) \arrow[swap]{rrr}{\alpha_J} \arrow[swap]{ur}{\Gamma(j,J)} &&& \Delta(J,J) \arrow[swap]{ur}{\Delta(j,J)} \end{tikzcd}
\]
\caption{The strong dinaturality condition (\cref{tag:EQstrongDinat}): if the diamond commutes, so too does the outer hexagon. }\label{tag:FIGstrongDinat}

\end{figure}

\mkDefn{DinaturalTransformation}

\begin{proposition}
Every strong dinatural transformation is a dinatural transformation.
\end{proposition}

\begin{note}
If the domain difunctor $\Gamma$ is either a presheaf $\C\op\to\SET$ padded with a dummy covariant argument, or a functor $\C\to\SET$ padded with a dummy contravariant argument, then every dinatural transformation $\Gamma\tod\Delta$ is a strong dinatural transformation $\Gamma\tos\Delta$.
\end{note}

\begin{definition} \label{tag:StructDefn}

Given $\Gamma\colon\C\op\times\C\to\SET$, define the category of $\Gamma$-structures, written $\Struct{\Gamma}$, to have objects $(I\colon\abs{\C},x\colon \Gamma(I,I))$ and morphisms $(I,x)$ to $(J,y)$ to be morphisms $j\colon\C\:[I,J]$ such that $$ \Gamma(I,j)\;x = \Gamma(j,J)\;y.$$

Write $\mc U_\Gamma\colon\Struct{\Gamma} \to \C$ for the first projection functor $(I,x)\mapsto I$.
\end{definition}

\begin{proposition} \label{tag:0b5f90}

The functor $\mc U_\Gamma$ is faithful.
\end{proposition}

\begin{proposition} \label{tag:b7b5ca}

The following data are equivalent:
\begin{itemize}
\item a strong dinatural transformation $\alpha\colon\Gamma\tos\Delta$
\item a functor $A\colon\Struct\Gamma\to\Struct\Delta$ such that $\mc U_\Delta\circ A = \mc U_\Gamma$.
\end{itemize}
\end{proposition}

\mkDefn{coe-difunctor}

\pdfImg{coe-difunct}{The $\coe[\Gamma]$ operation}{coe-difunct}

\begin{proposition} \label{tag:255510}

The following properties hold for any $\Gamma$, $x$, and any isos $j$, $k$
\begin{itemize}
\item $\coe[\Gamma](\iden,x) = x$
\item $\coe[\Gamma](k\circ j,x) = \coe[\Gamma](k,\coe[\Gamma](j,x))$
\item $j$ is a $\Struct\Gamma$-isomorphism $(I,x) \cong (J,\coe(j,x))$.
\end{itemize}
\end{proposition}

\begin{proposition} \label{tag:693064}

For any dinatural transformation $\alpha\colon\Gamma\tod\Delta$,\footnote{In particular: any \textit{strong} dinatural transformation $\alpha\colon\Gamma\tos\Delta$.} any $x\colon \Gamma(I,I)$ and any $j\colon I\cong J$, $$ \alpha_J\;(\coe[\Gamma](j,x)) = \coe[\Delta](j,\alpha_I\;x).$$
\end{proposition}

\begin{proof}
Let $$ w\dEqual \Gamma(j\inv,I)\;x \qquad \colon \Gamma(J,I). $$ Then $x=\Gamma(j,I)\;w$ and, by the dinaturality of $\alpha$, we know that $$ \Delta(I,j)\;(\alpha_I\;x) = \Delta(j,J)\;(\alpha_J\;(\Gamma(J,j)\;w)), $$ i.e. $$ \Delta(I,j)\;(\alpha_I\;x) = \Delta(j,J)\;(\alpha_J\;(\coe[\Gamma](j,x))). $$ Now apply $\Delta(j\inv,J)$ to both sides and get $$ \coe[\Delta](j,\alpha_I\;x) = \alpha_J\;(\coe[\Gamma](j,x)).$$
\end{proof}

\begin{figure}
\[
    \begin{tikzcd}[sep=large] & \Gamma(I,I) \arrow{rr}{\alpha_I}\arrow{dd}{\coe[\Gamma](j)}  \arrow[shift right=1,near end]{dl}{\Gamma(j\inv,I)} && \Delta(I,I) \arrow{dr}{\Delta(I,j)} \arrow[swap]{dd}{\coe[\Delta](j)} \\
        \Gamma(J,I) \arrow[shift left=3]{ur}{\Gamma(j,I)} \arrow[swap]{dr}{\Gamma(J,j)} && && \Delta(I,J) \arrow[swap,shift left=1,near start]{dl}{\Delta(j\inv,J)} \\
 & \Gamma(J,J) \arrow[swap]{rr}{\alpha_J} && \Delta(J,J) \arrow[swap,shift right=3]{ur}{\Delta(j,J)} \end{tikzcd}
\]
\caption{{}\cref{tag:693064} for $\alpha\colon\Gamma\tod\Delta$ and $j\colon I\cong J$.}\label{tag:FIGdinatCoe}

\end{figure}

\section{Generalized Difunctors and Paranaturality}\label{preliminaries}

\mkDefn{faithfulDisplayedCategory}

\mkDefn[\formalLink{GenDifunctor}]{GenDifunctor}

\begin{note}

Especially when the underlying category $\C$ is $\SET$, we'll think of the a generalized difunctor $\mb D$ as encoding a notion of ``structure'':
\begin{itemize}
\item $\mb D_I$ is the collection of ``structures based on $I$''
\item $(x,y)\in \mb D_j$ means that the function $j$ (from the underlying set of $x$ to the underlying set of $y$) is a ``structure-preserving morphism'' from $x$ to $y$.
\item $\coe[\mb D]$ says that we can rename the elements of the underlying sets without changing the $\mb D$-structure
\end{itemize}

\end{note}

\begin{example} \label{tag:difunctorsAreGenDifunctors}

Every difunctor $\Gamma$ over $\C$ can be viewed as a generalized difunctor:
\begin{itemize}
\item $\Gamma_I\dEqual \Gamma(I,I)$
\item $(x,y)\in \Gamma_j$ iff
  $$ \Gamma(I,j)\;x = \Gamma(j,J)\;y $$
\item $\coe[\Gamma]$ as defined in \defnRef{coe-difunctor}
\end{itemize}

\end{example}

\begin{example} \label{tag:monoidGenDifunctor}

Consider the difunctor $\mu\colon\SET\op\times\SET\to\SET$ with object part $$ \mu(X,Y) \dEqual Y \times (X\to X\to Y). $$ An object of $\Struct{\mu}$ is a set $X$ with a distinguished element $e\colon X$ and binary operation $m\colon X\to X\to X$. The category $\MON$ of \textit{monoids} is a subcategory of $\Struct{\mu}$---the full subcategory of those $\mu$-structures $(X,e,m)$ such that $m$ is associative and $e$ is a unit for $m$. However, there is (apparently) no way to to realize $\MON$ directly as $\Struct{\Gamma}$ for some $\Gamma$: the definition of a difunctor must separate the ``polarities'', as we have done in the definition of $\mu$ above ($X$ being the ``negatively''-positioned elements and $Y$ the ``positively''), hence equations like $m(x,e)=x$ become ill-typed.

However, $\MON$ \textit{is} a generalized difunctor: for any set $X$, $\MON_X$ is the set of monoid structures on $X$---pairs $(e\colon X, m\colon X \to X \to X$) where $m$ is associative and $e$ is a unit for $m$; $((e,m),(e',m'))\in \MON_f$ for $f\colon X\to Y$ just in case $f$ is a monoid homomorphism from $(X,e,m)$ to $(Y,e',m)$; and the $\coe$ operation is just renaming the elements of a monoid structure on a set $X$ via a bijection $X\cong Y$.\footnote{ Observe that these definitions of $\MON_f$ and $\coe[\MON]$ are the ones we'd get by viewing the family of sets $\set{\MON_X}_{X\colon\SET}$  as a subfamily of $\set{\mu(X,X)}_{X\colon\SET}$ (since $\MON_X \subseteq \mu(X,X)$ for all $X$), and restricting $\mu_f$ and $\coe[\mu]$ to this subfamily, inducing the correct generalized difunctor structure on $\set{\MON_X}_{X\colon\SET}$.  } Equivalently, we can observe that the forgetful functor $\MON \to \SET$ is faithful and a split isofibration.
\end{example}

\begin{example} \label{tag:topGenDifunctor}

The collection of \textit{topologies} $\TOP$ is a generalized difunctor over $\SET$:
\begin{itemize}
\item $\TOP_I$ is the collection of topologies on the set $I$
\item $(\tau_1,\tau_2)\in \TOP_j$ if $j$ is continuous with respect to the topologies $\tau_1,\tau_2$
\item for a bijection $j\colon I\cong J$ and topology $\tau$ on $I$, $\coe(j,\tau)$ is a topology on $J$, whose open sets are just the $j$-images of open sets on $I$ (in the topology $\tau$).
\end{itemize}
\end{example}

\begin{note}
We can distinguish several nested classes of generalized difunctors:
\begin{itemize}
\item the category of elements of a functor $\C\to\SET$ (or presheaf $\C\op\to\SET$), viewed as a faithful displayed category (via \cref{tag:difunctorsAreGenDifunctors})
\item difunctors $\Gamma\colon\C\op\times\C\to\SET$, viewed as generalized difunctors (\cref{tag:difunctorsAreGenDifunctors})
\item generalized difunctors (like $\MON$ from \cref{tag:monoidGenDifunctor}) which are sub-generalized difunctors of difunctors obtained by imposing equations
\item generalized difunctors (like $\TOP$ from \cref{tag:topGenDifunctor}) which are not of an ``algebraic character'', i.e. (apparently) not the class of algebras for any generalized algebraic theory
\end{itemize}
\end{note}

\mkDefn{TotalCat}

\begin{example}\
\begin{itemize}
\item When $\mb D$ is a difunctor $\Gamma$, then $\Sigma\mb D$ is $\Struct\Gamma$ from \cref{tag:StructDefn}
\begin{itemize}
\item if $\Gamma$ is a functor $\C\to\SET$ or $\C\op\to\SET$, then $\Sigma\mb D$ is its \textit{category of elements};
\item if $\Gamma$ is the difunctor $\C\:[T(-),-]$ for some endofunctor $T\colon\C\to\C$, then $\Sigma\mb D$ is the \textit{category of $T$-algebras};
\item if $\Gamma$ is the difunctor $\C\:[-,T(-)]$ for some endofunctor $T\colon\C\to\C$, then $\Sigma\mb D$ is the \textit{category of $T$-coalgebras};
\end{itemize}
\item $\Sigma\MON$ and $\Sigma\TOP$ are just the usual categories of monoids and topological spaces, respectively.
\end{itemize}
\end{example}

\mkDefn[\formalLink{ParanaturalTransformation}]{ParanaturalTransformation}

\pdfImg{paranat-mor}{ A paranatural transformation $\alpha$ preserves morphisms: if there's a morphism over $j$ from $x$ to $y$, then there's one from $\alpha_I\;x$ to $\alpha_J\;y$}{paranat-mor}

\begin{note}
We'll write $$ \int_{K\colon\C} \mb D_K \dD \mb E_K $$ for the type of paranatural transformations from $\mb D$ to $\mb E$.
\end{note}

\begin{proposition}
Given difunctors $\Gamma,\Delta$, the following data are equivalent
\begin{itemize}
\item a strong dinatural transformation $\Gamma\tos\Delta$
\item a paranatural transformation from $\Gamma$ to $\Delta$ (viewed as generalized difunctors), i.e. an element of $$\int_I \Gamma(I,I)\dD \Delta(I,I).$$
\end{itemize}
\end{proposition}

\section{DiYoneda and Exponentiation}\label{diYoneda}

\subsection{DiYoneda}

\begin{example}
For any category $\C$, we have the hom difunctor $\hom\colon\C\op\times\C\to\SET$. $\Struct{\hom}$ is the category of $\C$\textit{-endofunctions}: the objects of $\Struct{\hom}$ are objects $I$ of $\C$ equipped with a morphism $\C\:[I,I]$. A $\Struct{\hom}$-morphism from $(I,x)$ to $(J,y)$ is a morphism $j\colon\C\:[I,J]$ such that $j\circ x = y \circ j$ (this latter equation is the definition of when $(x,y)\in \hom[j]$). The $\coe$ operation is obtained by convolution by the isomorphism: given $x\colon \C\:[I,I]$ and $j\colon I\cong J$, obtain $$ \coe(j,x) \dEqual j \circ x \circ j\inv \qquad \colon \C\:[J,J].$$
We can obtain the \textit{generalized} difunctor $\Auto$ of \textit{automorphisms} (that is, endo-isomorphisms) as a full sub-generalized difunctor of $\hom$: $\Auto_I$ is the set of $\C$-\textit{iso}morphisms $I\cong I$. Since isomorphisms are closed under composition, the above definition of $\coe(j,x)$ defines an automorphism of $J$ when $x$ is an automorphism of $I$. The notion of \enquote{morphism} between $x\colon I\cong I$ and $y\colon J \cong J$ is also inherited from $\Struct\hom$: a morphism $j\colon\C\:[I,J]$ such that $j\circ x=y\circ j$; note that this equation can hold without $j$ itself being an isomorphism.
\end{example}

    \begin{lemma}\label{state:DiYonedaLemma}

For each generalized difunctor $\mb E$, there is an isomorphism $$\mb E_I\cong \int_K I \cong K \dD \mb E_K$$
paranatural in $I$.

    \end{lemma}

The right-hand side is given by the following definition (see \cref{multivariable} for a more general development):
\mkDefn{DiYonedaTranspose}

Now the proof of \cref{state:DiYonedaLemma}:
\begin{proof}
For an object $I$ of $\C$, define $\Phi_I\colon \mb E_I \cong \widetilde{\mb E}_I$ by $$ (\Phi_I\;x)_K\;(k\colon I\cong K) \dEqual \coe[\mb E](k,x) $$ and $$ \Phi_I\inv\;\varphi \dEqual \varphi_I\;\iden_I. $$
We must check:
\begin{itemize}
\item $\Phi_I\;x$ is indeed an element of $\widetilde{\mb E}_I$:  for $j\colon I\cong J$ and $k'\colon J\cong K$, $$ \coe[\mb E](k',(\Phi_I\;x)_J\;j) = \coe[\mb E](k',\coe[\mb E](j,x)) = \coe(k'\circ j,x) = (\Phi_I\;x)_K\;(k'\circ j).$$
\item $\Phi$ is a paranatural transformation: \cref{tag:8b11f4}
\item $\Phi\inv$ is a paranatural transformation: \cref{tag:ecd4a8}
\item $\Phi$ and $\Phi\inv$ are inverses: since $\coe[\mb E](\iden,x)=x$, we have $\Phi_I\inv\circ\Phi_I$ is the identity and, conversely
\begin{align*}
(\Phi_I(\Phi_I\inv\;\varphi))_K\; k
&= \coe[\mb E](k,\Phi_I\inv\;\varphi)\\
&= \coe[\mb E](k,\varphi_I\;\iden_I)\\
&= \varphi_K\;(k\circ\iden_I)\\
&= \varphi_K\;k
\end{align*}
so $\Phi_I(\Phi_I\inv\;\varphi) = \varphi$ too.
\end{itemize}
\end{proof}

\subsection{Exponentials}

\begin{example} \label{tag:prodDefn}

For generalized difunctors $\mb D,\mb E$, define the generalized difunctor $\mb D\times\mb E$
\begin{itemize}
\item $(\mb D\times\mb E)_I \dEqual \mb D_I\times\mb E_I$
\item $((d,e),(d',e'))\in(\mb D\times\mb E)_j$ iff $(d,d')\in\mb D_j$ and $(e,e')\in \mb E_j$
\item $\coe(j,(d,e)) \dEqual (\coe(j,d),\coe(j,e))$
\end{itemize}

\end{example}

\begin{proposition}
The projection maps $\pi_1\colon \mb D\times\mb E \to \mb D$ and $\pi_2\colon\mb D\times\mb E\to \mb E$ are paranatural transformations, and satisfy the universal mapping property of the product.
\end{proposition}

\begin{definition}
Given generalized difunctors $\mb D,\mb E$ over $\C$, define the generalized difunctor $\mb E^{\mb D}$
\begin{itemize}
\item an element $\varphi$ of $\mb E^{\mb D}_I$ is a family of functions $$ \varphi_K \colon I\cong K \to \mb D_K \to \mb E_K$$ for each object $K$ of $\C$ such that,
    \begin{itemize}
    \item for $j\colon I\cong J$ and $f\colon J\cong K$ and $d'\colon\mb D_J$,
        \begin{equation} \label{tag:EQexponPresCoe}
        \coe[\mb E](f,\varphi_J\;j\;d') = \varphi_K\;(f\circ j)\;(\coe[\mb D](f,d'))
        \end{equation}

    \end{itemize}
\item $(\varphi,\psi)$ is in $\mb E^{\mb D}_j$ iff for all $(d,d')\in \mb D_j$ $$ (\varphi_I\;\iden_I\;d,\psi_J\;\iden_J\;d') \in \mb E_j$$ (see \cref{tag:exponMorphiEquiv})

\item for $j\colon I\cong J$ and $\varphi\colon\mb E^{\mb D}_I$, define $\coe[{\mb E}^{\mb D}](j,\varphi)\colon{\mb E}^{\mb D}_J$ by
$$
(\coe[{\mb E}^{\mb D}](j,\varphi))_K\;(f\colon J\cong K)\;(dd\colon\mb D_K) \dEqual \varphi_K\;(f\circ j)\;dd
$$

\end{itemize}
\end{definition}

\begin{proposition} \label{tag:exponMorphiEquiv}

For $j\colon \C\:[I,J]$ and $\varphi\in\mb E^{\mb D}_I$, $\psi\in\mb E^{\mb D}_J$, the following are equivalent:
\begin{itemize}
\item for all $g\colon I\cong K$, $h\colon J\cong L$, and all $(dd,dd')\in \mb D_{h\circ j\circ g\inv}$ $$ (\varphi_K\;g\;dd,\psi_L\;h\;dd') \in \mb E_{h\circ j\circ g\inv} $$
\item for all $(d,d')\in \mb D_j$, $$(\varphi_I\;\iden_I\;d,\psi_J\;\iden_J\;d') \in \mb E_j$$
\end{itemize}

\end{proposition}
\begin{proof}
The second statement is an instance of the first ($g=\iden_I$, $h=\iden_J$). Suppose the second holds and let $g,h,dd,dd'$ be given. Then let $$ d\dEqual\coe[\mb D](g\inv,dd) \colon \mb D_I \qquad\text{and}\qquad d'\dEqual\coe[\mb D](h\inv,dd') \colon \mb D_J. $$
Now, since $(d,dd)\in\mb D_g$ and $(dd,dd')\in \mb D_{h\circ j\circ g\inv}$ and $(dd',d)\in\mb D_{h\inv}$, we have by the heterogeneous transitivity that $$ (d,d') \in \mb D_{h\inv\circ(h\circ j\circ g\inv)\circ g} \qquad\text{i.e.}\quad (d,d')\in\mb D_j.$$
Then, by hypothesis, this implies that $$ (\varphi_I\;\iden_I\;d,\psi_J\;\iden_J\;d') \in \mb E_j $$
i.e. $$ (\varphi_I\;(g\inv\circ g)\;(\coe[\mb D](g\inv,dd)),\psi_J\;(h\inv \circ h)\;(\coe[\mb D](h\inv,dd'))) \in \mb E_j $$
Now, by \cref{tag:EQexponPresCoe}, $$ \varphi_I\;(g\inv\circ g)\;(\coe[\mb D](g\inv,dd)) = \coe[\mb E](g\inv,\varphi_K\;g\;dd)$$ and likewise $\psi_J\;(h\inv\circ h)\;(\coe[\mb D](h\inv,dd')) = \coe[\mb E](h\inv, \psi_L\;h\;dd')$, so we can now write $$ (\coe[\mb E](g\inv,\varphi_K\;g\;dd,\quad \coe[\mb E](h\inv,\psi_L\;h\;dd')) \quad\in\quad \mb E_j$$
However, we know that
\begin{align*}
&(\varphi_K\;g\;dd,\quad\coe[\mb E](g\inv,\varphi_K\;g\;dd)) &\in \mb E_{g\inv}\\
&(\coe[\mb E](h\inv,\varphi_L\;h\;dd'),\quad \psi_L\;h\;dd') &\in \mb E_{h}
\end{align*}
so, combining these by heterogeneous transitivity, $$ (\varphi_K\;g\;dd,\psi_L\;h\;dd') \quad\in\quad \mb E_{h\circ j\circ g\inv} $$

\end{proof}

\begin{proposition} \label{tag:exponObj}

For any generalized difunctors $\mb D,\mb E$ and any object $I$, there is a bijection $$ \mb E^{\mb D}_I \quad\cong\quad \mb D_I \to \mb E_I.$$
\end{proposition}

\begin{proof}
Denote the forward direction $\_\$\_$, defined by: $$ \left(\varphi\colon\mb E^{\mb D}_I\right)\,\$\,(d\colon \mb D_I) \dEqual \varphi_I\;\iden_I\;d \quad\colon\mb E_I $$
and the other direction $\Lambda$: $$ (\Lambda\; (p\colon \mb D_I\to \mb E_I))_K \;(g\colon I\cong K)\;(dd\colon \mb D_K) \dEqual \coe[\mb E](g, p(\coe[\mb D](g\inv,dd))).$$
That $(\Lambda\;p)\,\$\,d = p\;d$ for all $p,d$ is immediate. Conversely, we have $\Lambda (\varphi\,\$\,\_) = \varphi$, where $(\varphi\,\$\,\_)$ is the function $(d\colon \mb D_I \mapsto \varphi\,\$\,d) \quad\colon \mb D_I \to \mb E_I$: for arbitrary $K,g,dd$,
\begin{align*}
(\Lambda(\varphi\,\$\,\_))_K\;g\;dd
&= \coe[\mb E](g,\;\varphi\,\$\,(\coe[\mb D](g\inv,dd))) \tag{Defn. $\Lambda$}\\
&= \coe[\mb E](g,\;\varphi_I\;\iden_I\;(\coe[\mb D](g\inv,dd))) \tag{Defn. $\$$}\\
&= \coe[\mb E](g,\;\varphi_I\;(g\inv\circ g)\;(\coe[\mb D](g\inv,dd)))\\
&= \coe[\mb E](g,\coe[\mb E](g\inv,\varphi_K\;g\;dd))
\tag{{}\cref{tag:EQexponPresCoe}{}}\\
&= \varphi_K\;g\;dd.
\end{align*}

\end{proof}

\begin{proposition} \label{tag:exponMorph}

For any $p\colon\mb D_I \to \mb E_I$ and $q\colon\mb D_J\to\mb E_J$ and $j\colon\C\:[I,J]$
\[  (\Lambda\;p,\Lambda\;q) \in \mb E^{\mb D}_j \qquad\text{iff}\qquad \forall (d,d')\in \mb D_j,\;(p\;d,q\;d')\in\mb E_j. \]

\end{proposition}

\begin{note}
Suppose $\alpha$ is a paranatural transformation whose domain is of the form $\mb E^{\mb D}$, i.e.
\[  \alpha \colon \int_K \mb E^{\mb D}_K \dD \mb F_K \]
for some generalized difunctor $\mb F$. \textit{Of what}, we might ask, does this $\alpha$ consist?

The underlying family of maps for $\alpha$ is of shape $$ \alpha_I \colon\quad \left(\varphi\colon \mb E^{\mb D}_I\right)\quad \mapsto \quad (\alpha_I\;\varphi \colon \mb F_I).$$ \cref{tag:exponObj} says that every such $\varphi$ is of the form $\Lambda p$ for some $p\colon \mb D_I \to \mb E_I$.

Next, we can slightly unpack the requirement that $\alpha$ preserves morphisms (if $(\Lambda p,\Lambda q)\in\mb E^{\mb D}_j$, then $(\alpha_I(\Lambda p),\alpha_J(\Lambda q))\in \mb F_j$): recall \cref{tag:exponMorph}, which tells us that $(\Lambda p,\Lambda q)\in\mb E^{\mb D}_j$ means that, for all $(d,d')\in \mb D_j$, $(p\;d,q\;d')\in \mb E_j$. So, ``$\alpha$ preserves morphisms'' means that if $p,q$ behave like the parallel components of a paranatural transformation $\mb D\to\mb E$ with respect to $j$---that is, if they always send pairs $(d,d')$ of $\mb D$-structures between whom $j$ is a $\mb D$-structure-preserving map ($(d,d')\in\mb D_j$) to pairs $(p\;d,q\;d)$ between whom $j$ is a $\mb E$-structure-preserving map---then we can conclude $j$ is a $\mb F$-structure-preserving map between $\alpha_I\;(\Lambda p)$ and $\alpha_J(\Lambda q)$.

Finally, to say that $\alpha$ preserves $\coe$ is the statement: for $j\colon I\cong J$, $$\alpha_J(\coe[\mb E^{\mb D}](j,\Lambda p))=\coe[\mb F](j,\alpha_I(\Lambda p)).$$ Observe that $\coe[\mb E^{\mb D}](j,\Lambda p) = \Lambda (\coe[\mb E](j) \circ p \circ \coe[\mb D](j\inv))$:

\begin{align*}
&(\coe[\mb E^{\mb D}](j,\Lambda p))_K\;(f\colon J\cong K)\;(dd\colon\mb D_K)\\
&\dEqual (\Lambda p)_K\;(f\circ j)\;dd\\
&= \coe[\mb E](f\circ j, p(\coe[\mb D](j\inv\circ f\inv, dd)))\\
&= \coe[\mb E](f,\coe[\mb E](j, p(\coe[\mb D](j\inv,\coe[\mb D](f\inv, dd)))))\\
&= \coe[\mb E](f,\quad \big(\coe[\mb E](j) \circ p \circ \coe[\mb D](j\inv)\big)\;(\coe[\mb D](f\inv, dd)) )\\
&= (\Lambda \big(\coe[\mb E](j) \circ p \circ \coe[\mb D](j\inv)\big))_K\;f\;dd
\end{align*}
\end{note}

\begin{example} \label{tag:trickySecondOrder}

Consider the difunctors
\begin{itemize}
\item $\hom\colon\SET\op\times\SET\to\SET$, the hom-set difunctor;
\item $\Iden\colon\SET\to\SET$, the identity functor, padded with a dummy variable to be a difunctor $\SET\op\times\SET\to\SET$
\end{itemize}
and consider them as generalized difunctors over $\SET$. Then apply \cref{tag:exponObj} and \cref{tag:exponMorph}:
\[  \left(\Iden^{\hom}\right)_I \qquad\cong\qquad (I\to I)\to I  \]
and, for $p\colon(I\to I)\to I$ and $q\colon(J\to J)\to J$ and $j\colon I\to J$,
\[ (\Lambda\;p,\Lambda\;q) \in \mb E^{\mb D}_j \qquad\text{iff}\qquad \forall x\colon I\to I, y\colon J\to J,\; j\circ x = y \circ j \implies j(p\;x)= q\;y. \]
\end{example}

\begin{example}
Consider ${\hom}^{\hom}$. Again, by \cref{tag:exponObj},
\[  \left({\hom}^{\hom}\right)_I \qquad\cong\qquad (I \to I) \to (I \to I). \]
Given a global element $\alpha\colon\smallint_I \UNIT \dD ({\hom}^{\hom})_I$, we can associate a natural number $$ N(\alpha) \dEqual ((\alpha_{\N}\;*)\,\$\,\mathsf{succ})\;0 $$ (where $*$ is the unique element of $\UNIT$) and, conversely, for any $n\colon\N$, we can form $\overline n\colon\smallint_I \UNIT\dD({\hom}^{\hom})_I$ by $$ (\overline n)_I\;* \dEqual \Lambda (\lambda i \to i^n). $$ It's immediate to see that $N(\overline n)=n$.
\end{example}

\begin{proposition}
There is a natural isomorphism
\[  \int_K \mb B_K \times \mb D_K \dD \mb E_K \quad\cong\quad \int_K \mb B_K \dD \mb E^{\mb D}_K.  \]
\end{proposition}

\section{The Structural End Calculus}\label{structEnd}

\begin{proposition}
For any generalized difunctor $\mb D$ and difunctor $\Delta$ (over the same category $\C$), $$ \int_{I\colon\C} \mb D_I\dD \Delta(I,I) \qquad= \qquad \int_{\mc D\colon\Sigma\mb D}\Delta(\mc U_{\mb D}\mc D,\mc U_{\mb D}\mc D). $$
\end{proposition}
The right-hand side is an \textit{end} in the sense of \cite[Chapter IX]{CWM}, namely the end of the difunctor $\Delta\colon\C\op\times\C\to\SET$ pre-composed in both arguments by the functor $\mc U_{\mb D}\colon\Sigma\mb D\to\C$ (see \defnRef{TotalCat}). A \textit{wedge}~\cite[Defn. 1.1.4]{loregian2023coend} for this functor consists of a set $V$ and a family of maps $$ \gamma_{(J,y)} \colon V \to \Delta(J,J) $$indexed over objects $(J,y)$ of $\Sigma\mb D$ which is appropriately dinatural: for any $j\colon\Sigma\mb D\:[(I,x),(J,y)]$ and any $v\colon V$, $$ \Delta(I,j)\;(\gamma_{(I,x)}\;v) = \Delta(j,J)\;(\gamma_{(J,y)}\;v).$$
A morphism of wedges from $(V,\gamma)$ to $(W,\epsilon)$ is a function $u\colon V\to W$ so that $\gamma_{(J,y)} = \epsilon_{(J,y)}\circ u$ for every $(J,y)$. The claim above---that $\smallint_I\mb D_I\dD \Delta(I,I)$ ``is'' the end of $\Delta$ pre-composed with $\mc U_{\mb D}$---is that $\smallint_I\mb D_I\dD\Delta(I,I)$ is the terminal wedge.
\begin{proof}
For convenience, write $X$ for the set $\smallint_I\mb D_I\dD\Delta(I,I)$ of all paranatural transforms from $\mb D$ to $\Delta$. First, we must make this a wedge $\omega_{(J,y)}\colon X\to \Delta(J,J)$. The definition must be: $$ \omega_{(J,y)}\;\alpha \dEqual \alpha_J\;y. $$ The dinaturality of this choice is actually the \textit{para}naturality of $\alpha$ (specifically that it must preserve morphisms): saying $j\colon\Sigma\mb D\:[(I,x),(J,y)]$ just means $(x,y)\in\mb D_j$, so we get that $(\alpha_I\;x,\alpha_J\;y)\in\Delta_j$, which is precisely the same as saying $$ \Delta(I,j)\;(\omega_{(I,x)}\;\alpha) = \Delta(j,J)\;(\gamma_{(J,y)}\;\alpha).$$
So we have a wedge. Now suppose $(V,\gamma)$ is some arbitrary wedge and pick arbitrary $v\colon V$. Then we define $\chi(v)\colon\smallint_I\mb D_I\dD \Delta(I,I)$ by $$ \chi(v)_I\;(x\colon\mb D_I) \dEqual \gamma_{(I,x)}\;v \qquad\colon\Delta(I,I). $$ To see that $\chi(v)$ preserves morphisms, pick $j$ and $(x,y)\in\mb D_j$ and notice that, since $(V,\gamma)$ is a wedge, we must have $(\gamma_{(I,x)}\;v, \gamma_{(J,y)}\;v) \in \Delta_j$. For preservation of $\coe[\mb D]$, let $j$ be an iso. Then, since $(x,\coe[\mb D](j,x))\in\mb D_j$, we again use the wedge condition for $\gamma$ to get $$ \Delta(I,j)\;(\gamma_{(I,x)}\;v) = \Delta(j,J)\;(\gamma_{(J,\coe[\mb D](j,x))}\;v)$$ and then apply $\Delta(j\inv,J)$ to either side and get
$$ \coe[\Delta](j,\chi(v)_I\;x) = \chi(v)_J\;(\coe[\mb D](j,x))$$ as desired. So we have a map $\chi\colon V\to X$. To see that $\chi$ is a wedge morphism, observe $$ (\omega_{(J,y)} \circ \chi)\;v = \chi(v)_J\;y = \gamma_{(J,y)}\;v$$ by definition. Finally, $\chi$ is the unique wedge morphism $(V,\gamma)\to (X,\omega)$: if $\eta\colon V\to X$ were such that $\omega_{(J,y)}\circ \eta = \gamma_{(J,y)}$ for all $(J,y)$, then, for any $v\colon V$ and any $(I,x)$, $$ \eta(v)_I\;x =: \omega_{(I,x)}\;(\eta(v)) = (\omega_{(I,x)}\circ \eta)\;v = \gamma_{(I,x)}\;v =: \chi(v)_I\;x$$ so $\eta = \chi$.
\end{proof}

The following generalizes \cite[Proposition 1]{uustalu2010} slightly.
\mkDefn{InitialAlgebraGenDifunctor}
\begin{note}
An initial algebra for $\mb D$ is indeed an initial object in the category $\Sigma\mb D$ from \defnRef{TotalCat}. In particular, if $\mb D$ is the difunctor $\C\:[T(-),-]$ for some endofunctor $T$, then an initial algebra for $\mb D$ is precisely an initial $T$-algebra in the usual sense.
\end{note}

    \begin{lemma}\label{state:UustaluYonedaLemma}

For any generalized difunctor $\mb D$ over $\C$ with initial algebra $(\mu_{\mb D},\inn_{\mb D})$, and any functor $K\colon\mb C\to\SET$, there is an isomorphism $$ \int_I \mb D_I\dD K(I) \cong K(\mu_{\mb D}).$$

    \end{lemma}

    \begin{proof}
The left-to-right direction sends $\alpha$ to $\alpha_{\mu_{\mb D}}\;\inn_{\mb D}$, and the opposite direction sends $z\colon K(\mu_{\mb D})$ to the paranatural transformation $\alpha^z$ given by $$ \alpha^z_I\;(x\colon\mb D_I) \dEqual K\;(\rec\;x)\;z. $$
This is indeed paranatural: if $(x,y)\in\mb D_j$ for some morphism $j$, then the requirement $(\alpha^z_I\;x,\alpha^z_J\;y)\in K_j$ is just the equation $$ K\;j\;(K\;(\rec\;x)\;z) = K\;(\rec\;y)\;z $$ which follows from the uniqueness of $\rec\;y$: since $(\inn_{\mb D},x)\in\mb D_{\rec\;x}$ and $(x,y)\in\mb D_j$, conclude $(\inn_{\mb D},y)\in \mb D_{j\circ\rec\;x}$, so it must be the case that $j\circ\rec\;x = \rec\;y$. For preservation of $\coe$, observe that $\coe[K]$ is just the operation sending $j\colon I\cong J$ and $w\colon K(I)$ to $K\;j\;w\colon K(J)$, so the $\coe$-preservation requirement for $\alpha^z$ is just that $$ K\;(\rec(\coe[\mb D](j,x)))\;z = K\;j\;(K\;(\rec\;x)\;z) $$ but $\rec(\coe[\mb D](j,x)) = j\circ\rec\;x$ by the same logic as before (and the requirement that $(x,\coe[\mb D](j,x))\in\mb D_j$).

Given $z\colon K(\mu_{\mb D})$, notice that $(\alpha^z)_{\mu_{\mb D}}\;\inn_{\mb D} = K\;(\rec\;\inn_{\mb D})\;z$, which must be $z$ by the (easily established) fact that $\rec\;\inn_{\mb D}$ must be the identity morphism on $\mu_{\mb D}$. So we have verified one composition gives the identity. Conversely, given $\alpha\colon\smallint_I \mb D_I\dD K(I)$ and arbitrary $I,x$, we can deduce $$ K\;(\rec\;x)\;(\alpha_{\mu_{\mb D}}\;\inn_{\mb D}) = \alpha_I\;x $$ by the fact that $\alpha$ preserves morphisms: since $(\inn_{\mb D},x)\in \mb D_{\rec\;x}$ by initiality, it follows that $(\alpha_{\mu_{\mb D}}\;\inn_{\mb D}, \alpha_I\;x) \in K_{\rec\;x}$, which is precisely the previous equation.

    \end{proof}

The instances contemplated by Uustalu were those where $\mb D$ was the difunctor $\C\:[T(-),-]$ for some endofunctor $T$. But the statement given here applies even when $\mb D$ is not a difunctor, such as the following.
    \begin{corollary}\label{state:monoidChoice}
There is only one operation $\alpha$ which, for each monoid $M$, picks out an element $\alpha_M$ of $M$ such that $f(\alpha_M)=\alpha_N$ for every monoid homomorphism $f\colon M\to N$.

    \end{corollary}

    \begin{corollary}\label{state:ringChoice}
There is a bijection between
\begin{itemize}
\item the collection of operations $\alpha$ which pick an element $\alpha_R$ of each ring $R$ in such a way that $f(\alpha_R)=\alpha_S$ for each ring homomorphism $f\colon R\to S$; and
\item the set $\mb Z$ of integers.
\end{itemize}

    \end{corollary}

\bibliographystyle{alpha}
\bibliography{biblio}

\appendix

\section{Incorrect Claims in Version 1}\label{errata}

\subsection{Lemma 2.13 (DiYoneda)}
The ``DiYoneda Lemma'' that appeared in Version 1~\cite[Defn. 2.13]{neumann2023paranatural} claimed that, for any difunctor $\Gamma$, there was a strong dinatural isomorphism
\begin{equation} \label{tag:originalDiYoneda}
\Gamma(I,J) \cong \int_K \C\:[J,K]\times \C\:[K,I] \dD \Gamma(K,K).
\end{equation}
The ``proof'' correctly notes that $x\colon\Gamma(I,I)$ can be sent to the strong dinatural transformation $\alpha(x)$, defined by $$ \alpha(x)_K\;(g\colon\C\:[I,K], h\colon\C\:[K,I]) \dEqual \Gamma(h,g)\;x. $$
This operation can be undone by the operation sending $\phi$ to $\phi_I(\iden_I,\iden_I)\colon\Gamma(I,I)$, because $\alpha(x)_I(\iden,\iden) = \Gamma(\iden,\iden)\;x = x$. However, the converse is not true: starting with a $\phi$ on the right-hand side, there's no guarantee that $\alpha(\phi_I(\iden,\iden))$ is the same strong dinatural transformation as $\phi$. Indeed, \cite{neumann2024updates} presents a concrete counterexample: if $\Gamma$ is the $\hom$ difunctor on $\SET$ and $\phi$ the strong dinatural transformation sending $(g\colon\C\:[I,K], h\colon\C\:[K,I])$ to $h\circ g\circ h\circ g$, then $\phi_I(\iden,\iden)$ will be $\iden$ and $\alpha(\iden) \neq \phi$. Notice in the corrected DiYoneda Lemma, \cref{state:DiYonedaLemma}, that $g$ and $h$ are restricted to be inverses of each other, eliminating this possibility.

\subsection{Parametricity is Strong Dinaturality}

Section 3 of Version 1 asserts an equivalence between parametricity (specifically the free theorems of \cite{freeTheorems}). This is not so. The simplest counterexample (to our knowledge) is the System F type $$ \forall X. ((X\to X) \to X) \to X. $$
The ``free theorem'' obtained for this type by application of Wadler's method is this: for any $\phi$ of this type, and for any $j\colon I\to J$, any $p\colon (I\to I)\to I$ and $q\colon (J\to J)\to J$, $$ \big( \forall x\;y, j\circ x = y\circ j \quad\text{implies}\quad j(p\;x)=q\;y\big) \qquad\text{implies}\qquad j(\phi_I\;p) = \phi_J\;q.$$
On the other hand, $\phi$ is a strong dinatural transformation $\int_I ((I\to I)\to I) \dD I$ if, for all $j,p,q$ as above, $$ \big(\forall r\colon J\to I,\;j(p(r\circ j)) = q(j\circ r)\big) \qquad\text{implies}\qquad j(\phi_I\;p) = \phi_J\;q.$$
Notice that the conclusions match, but, since we're using the twisted exponential to interpret $(I\to I)\to I$ as a difunctor in $I$, the strong dinaturality condition obtains this unfortunate dependence on functions of type $J\to I$.

Note that \cref{tag:trickySecondOrder} uses our paranatural exponential for $\Iden^{\hom}$, and manages to get this right: $j\colon I\to J$ is a morphism from $p\colon(I\to I)\to I$ to $q\colon(J\to J)\to J$ just in case $j\circ x = y \circ j$ implies $j(p\;x)= q\;y$ for all $x,y$.

\section{Multi-variable Generalized Difunctors}\label{multivariable}

\begin{definition} \label{tag:BinaryGenDifunctDefn}

A \textbf{binary generalized difunctor} $\mb X$ over $\C$ is
\begin{itemize}
\item for each $I,J$, a set $\mb X_{I,J}$
\item for every $j\colon\C\:[I,J]$ and $\ell\colon\C\:[K,L]$, a relation $$ \mb X_{j,\ell} \subseteq \mb X_{I,K} \times \mb X_{J,L}$$ which is a heterogeneous preorder in each argument
\item operators
\begin{align*}
\coe[\mb X-,K] &\colon (I \cong J) \to \mb X_{I,K} \to \mb X_{J,K}\\
\coe[\mb XI,-] &\colon (K \cong L) \to \mb X_{I,K} \to \mb X_{I,L}
\end{align*}
each satisfying the usual rules for $\coe$.
\end{itemize}

\end{definition}

\section{Faithful Isofibrations}\label{fibration}

\begin{definition} \label{tag:splitIsofibDefn}

A  \textbf{split isofibration} is a functor $\mc F\colon\mb E\to\mb C$ along with
\begin{itemize}
\item for each object $e$ of $\mb E$ and $\mb C$-iso $j\colon\mc F(e)\cong J$ a given object $\coe(j,e)$ of $\mb E$ where $\mc F(\coe(j,e)) = J$
\item for each $e$ and $j$, an iso $\coh(j,e)\colon e \cong \coe(j,e)$ where $\mc F(\coh(j,e))= j$
\end{itemize}
such that
\begin{itemize}
\item $\coe(\iden,e) = e$ and $\coe(k\circ j,e)=\coe(k,\coe(j,e))$
\item $\coh(\iden,e)=\iden$ and $\coh(k\circ j,e)=\coh(k,\coe(j,e))\circ\coh(j,e)$.
\end{itemize}

\end{definition}

\begin{proposition} \label{tag:2c9bfc}

The following data are equivalent.
\begin{itemize}
\item A faithful, split isofibration over $\C$
\item A generalized difunctor over $\C$
\end{itemize}

\end{proposition}

\begin{proof}
The faithful split isofibration corresponding to a generalized difunctor $\mb D$ consists of the category $\Sigma\mb D$ defined in \defnRef{TotalCat}, along with its forgetful functor $\mc U_{\mb D}\colon\Sigma\mb D\to\C$. Conversely, the generalized difunctor $\mb F$ corresponding to a faithful split isofibration $\mc F$ is given by $\mb F_I\dEqual \mc F\inv(I)$, and the definition carries over in the obvious way.
\end{proof}

\begin{proposition}
The following data are equivalent
\begin{itemize}
\item A paranatural transformation from $\mb D$ to $\mb E$
\item A functor $A : \Sigma\mb D\to\Sigma\mb E$ such that
\begin{itemize}
\item $\mc U_{\mb E}\circ A = \mc U_{\mb D}$
\item $A(J,\coe[\mb D](j,d))=(J,\coe[\mb E](j,A(I,d))$ for all $j,d$
\end{itemize}
\end{itemize}
\end{proposition}

\section{Subproofs}\label{subproofs}

\subsection{DiYoneda}

\begin{lemma} \label{tag:8b11f4}

The assignment $\Phi_I \colon x \mapsto (\lambda\;K\;k \to \coe[\mb E](k,x))$ is a paranatural transformation $\smallint_K \mb E_K \dD \widetilde{\mb E}_K$

\end{lemma}
\begin{proof}
\begin{itemize}
\item if $j\colon\C\:[I,J]$ and $(x,y)\in\mb E_j$, then, since $$ (\Phi_I\;x)_I\;\iden_I \dEqual \coe[\mb E](\iden_I,x) = x $$ and likewise $(\Phi_J\;y)_J\;\iden_J = y$, we have $$ ((\Phi_I\;x)_I\;\iden_I, (\Phi_J\;y)_J\;\iden_J) \in \mb E_j $$ and thus $(\Phi_I\;x,\Phi_J\;y)\in\widetilde{\mb E}_j$;
 \item if $j\colon I\cong J$ and $x\colon\mb E_I$, then, for any $k\colon J\cong K$,
\begin{align*}
(\Phi_J\;(\coe[\mb E](j,x)))_K\;k
&\dEqual \coe[\mb E](k,\coe[\mb E](j,x))\\
&= \coe[\mb E](k\circ j,x)\\
&= (\Phi_I\;x)_K\;(k\circ j)\\
&= (\coe[\sim](j,\Phi_I\;x))_K\;k
\end{align*}
so $\Phi_J\;(\coe[\mb E](j,x)) = \coe[\sim](j,\Phi_I\;x)$.
\end{itemize}
\end{proof}

\begin{lemma}  \label{tag:ecd4a8}

The assignment $\Phi\inv_I \colon \varphi \mapsto \varphi_I\;\iden_I$ is a paranatural transformation $\smallint_K \widetilde{\mb E}_K \dD \mb E_K$

\end{lemma}
\begin{proof}
\begin{itemize}
\item if $(\varphi,\psi)\in\widetilde{\mb E}_j$, then by definition of $\widetilde{\mb E}$, we have $(\varphi_I\;\iden_I,\psi_J\;\iden_J)\in\mb E_j$, i.e. $(\Phi_I\inv\;\varphi,\Phi_J\inv\;\psi)\in\mb E_j$;
\item if $j\colon I\cong J$ and $\varphi\colon\widetilde{\mb E}_I$,
\begin{align*}
\Phi_J\inv\;(\coe[\sim](j,\varphi))
&= (\coe[\sim](j,\varphi))_J\;\iden_J\\
&= \varphi_J\;(\iden_J\circ j) \tag{defn. of $\coe[\sim]$}\\
&= \coe[\mb E](j,\varphi_I\;\iden_I) \tag{$\varphi\colon\widetilde{\mb E}_I$}\\
&= \coe[\mb E](j,\Phi_I\inv\;\varphi)
\end{align*}
\end{itemize}
\end{proof}

\subsection{Exponentials}

\clearpage
\section{Agda Formalization}
\subsection{GenDifunctor.Core}

\begin{code}[hide]%
\>[0]\AgdaSymbol{\{-\#}\AgdaSpace{}%
\AgdaKeyword{OPTIONS}\AgdaSpace{}%
\AgdaPragma{--without-K}\AgdaSpace{}%
\AgdaPragma{--safe}\AgdaSpace{}%
\AgdaSymbol{\#-\}}\<%
\\
\\[\AgdaEmptyExtraSkip]%
\>[0]\AgdaKeyword{module}\AgdaSpace{}%
\AgdaModule{GenDifunctor.Core}\AgdaSpace{}%
\AgdaKeyword{where}\<%
\\
\\[\AgdaEmptyExtraSkip]%
\>[0]\AgdaKeyword{open}\AgdaSpace{}%
\AgdaKeyword{import}\AgdaSpace{}%
\AgdaModule{Level}\<%
\\
\>[0]\AgdaKeyword{open}\AgdaSpace{}%
\AgdaKeyword{import}\AgdaSpace{}%
\AgdaModule{Categories.Category}\<%
\\
\>[0]\AgdaKeyword{open}\AgdaSpace{}%
\AgdaKeyword{import}\AgdaSpace{}%
\AgdaModule{Relation.Binary.PropositionalEquality}\AgdaSpace{}%
\AgdaSymbol{as}\AgdaSpace{}%
\AgdaModule{≡}\AgdaSpace{}%
\AgdaKeyword{using}\AgdaSpace{}%
\AgdaSymbol{(}\AgdaOperator{\AgdaDatatype{\AgdaUnderscore{}≡\AgdaUnderscore{}}}\AgdaSymbol{)}\<%
\\
\>[0]\AgdaKeyword{open}\AgdaSpace{}%
\AgdaKeyword{import}\AgdaSpace{}%
\AgdaModule{Relation.Binary}\AgdaSpace{}%
\AgdaKeyword{using}\AgdaSpace{}%
\AgdaSymbol{(}\AgdaFunction{Rel}\AgdaSymbol{;}\AgdaSpace{}%
\AgdaFunction{REL}\AgdaSymbol{;}\AgdaSpace{}%
\AgdaRecord{IsEquivalence}\AgdaSymbol{)}\<%
\\
\>[0]\AgdaKeyword{import}\AgdaSpace{}%
\AgdaModule{Categories.Morphism}\<%
\\
\\[\AgdaEmptyExtraSkip]%
\>[0]\AgdaKeyword{private}\<%
\\
\>[0][@{}l@{\AgdaIndent{0}}]%
\>[2]\AgdaKeyword{variable}\<%
\\
\>[2][@{}l@{\AgdaIndent{0}}]%
\>[4]\AgdaGeneralizable{o}\AgdaSpace{}%
\AgdaGeneralizable{ℓ}\AgdaSpace{}%
\AgdaGeneralizable{e}\AgdaSpace{}%
\AgdaGeneralizable{o′}\AgdaSpace{}%
\AgdaGeneralizable{ℓ′}\AgdaSpace{}%
\AgdaSymbol{:}\AgdaSpace{}%
\AgdaPostulate{Level}\<%
\end{code}
\begin{formal}[\cref{defn:faithful displayed category}, \cref{defn:GenDifunctor}]\label{agda:GenDifunctor}\ 

\begin{code}%
\>[0]\AgdaKeyword{record}\AgdaSpace{}%
\AgdaRecord{GenDifunctor}\AgdaSpace{}%
\AgdaSymbol{(}\AgdaBound{C}\AgdaSpace{}%
\AgdaSymbol{:}\AgdaSpace{}%
\AgdaRecord{Category}\AgdaSpace{}%
\AgdaGeneralizable{o}\AgdaSpace{}%
\AgdaGeneralizable{ℓ}\AgdaSpace{}%
\AgdaGeneralizable{e}\AgdaSymbol{)}\AgdaSpace{}%
\AgdaSymbol{(}\AgdaBound{o′}\AgdaSpace{}%
\AgdaBound{ℓ′}\AgdaSpace{}%
\AgdaBound{e′}\AgdaSpace{}%
\AgdaSymbol{:}\AgdaSpace{}%
\AgdaPostulate{Level}\AgdaSymbol{)}\<%
\\
\>[0][@{}l@{\AgdaIndent{0}}]%
\>[2]\AgdaSymbol{:}\AgdaSpace{}%
\AgdaPrimitive{Set}\AgdaSpace{}%
\AgdaSymbol{(}\AgdaBound{o}\AgdaSpace{}%
\AgdaOperator{\AgdaPrimitive{⊔}}\AgdaSpace{}%
\AgdaBound{ℓ}\AgdaSpace{}%
\AgdaOperator{\AgdaPrimitive{⊔}}\AgdaSpace{}%
\AgdaBound{e}\AgdaSpace{}%
\AgdaOperator{\AgdaPrimitive{⊔}}\AgdaSpace{}%
\AgdaPrimitive{suc}\AgdaSpace{}%
\AgdaBound{o′}\AgdaSpace{}%
\AgdaOperator{\AgdaPrimitive{⊔}}\AgdaSpace{}%
\AgdaPrimitive{suc}\AgdaSpace{}%
\AgdaBound{ℓ′}\AgdaSpace{}%
\AgdaOperator{\AgdaPrimitive{⊔}}\AgdaSpace{}%
\AgdaPrimitive{suc}\AgdaSpace{}%
\AgdaBound{e′}\AgdaSymbol{)}\AgdaSpace{}%
\AgdaKeyword{where}\<%
\\
\>[2]\AgdaKeyword{eta-equality}\<%
\\
\>[2]\AgdaKeyword{private}\AgdaSpace{}%
\AgdaKeyword{module}\AgdaSpace{}%
\AgdaModule{C}\AgdaSpace{}%
\AgdaSymbol{=}\AgdaSpace{}%
\AgdaModule{Category}\AgdaSpace{}%
\AgdaBound{C}\<%
\\
\>[2]\AgdaKeyword{open}\AgdaSpace{}%
\AgdaModule{Categories.Morphism}\AgdaSpace{}%
\AgdaBound{C}\<%
\\
\\[\AgdaEmptyExtraSkip]%
\>[2]\AgdaKeyword{field}\<%
\\
\>[2][@{}l@{\AgdaIndent{0}}]%
\>[4]\AgdaField{Objᴰ}\AgdaSpace{}%
\AgdaSymbol{:}\AgdaSpace{}%
\AgdaFunction{C.Obj}\AgdaSpace{}%
\AgdaSymbol{→}\AgdaSpace{}%
\AgdaPrimitive{Set}\AgdaSpace{}%
\AgdaBound{o′}\<%
\\
\>[4]\AgdaField{Homᴰ}\AgdaSpace{}%
\AgdaSymbol{:}\AgdaSpace{}%
\AgdaSymbol{∀}\AgdaSpace{}%
\AgdaSymbol{\{}\AgdaBound{I}\AgdaSpace{}%
\AgdaBound{J}\AgdaSymbol{\}}\AgdaSpace{}%
\AgdaSymbol{→}\AgdaSpace{}%
\AgdaBound{C}\AgdaSpace{}%
\AgdaOperator{\AgdaFunction{[}}\AgdaSpace{}%
\AgdaBound{I}\AgdaSpace{}%
\AgdaOperator{\AgdaFunction{,}}\AgdaSpace{}%
\AgdaBound{J}\AgdaSpace{}%
\AgdaOperator{\AgdaFunction{]}}\AgdaSpace{}%
\AgdaSymbol{→}\AgdaSpace{}%
\AgdaFunction{REL}\AgdaSpace{}%
\AgdaSymbol{(}\AgdaField{Objᴰ}\AgdaSpace{}%
\AgdaBound{I}\AgdaSymbol{)}\AgdaSpace{}%
\AgdaSymbol{(}\AgdaField{Objᴰ}\AgdaSpace{}%
\AgdaBound{J}\AgdaSymbol{)}\AgdaSpace{}%
\AgdaBound{ℓ′}\<%
\\
\>[4]\AgdaOperator{\AgdaField{\AgdaUnderscore{}≈ᴰ\AgdaUnderscore{}}}\AgdaSpace{}%
\AgdaSymbol{:}\AgdaSpace{}%
\AgdaSymbol{∀}\AgdaSpace{}%
\AgdaSymbol{\{}\AgdaBound{I}\AgdaSymbol{\}}\AgdaSpace{}%
\AgdaSymbol{→}\AgdaSpace{}%
\AgdaFunction{Rel}\AgdaSpace{}%
\AgdaSymbol{(}\AgdaField{Objᴰ}\AgdaSpace{}%
\AgdaBound{I}\AgdaSymbol{)}\AgdaSpace{}%
\AgdaBound{e′}\<%
\\
\\[\AgdaEmptyExtraSkip]%
\>[4]\AgdaField{identity}\AgdaSpace{}%
\AgdaSymbol{:}\AgdaSpace{}%
\AgdaSymbol{∀}\AgdaSpace{}%
\AgdaSymbol{\{}\AgdaBound{I}\AgdaSymbol{\}}\AgdaSpace{}%
\AgdaSymbol{\{}\AgdaBound{x}\AgdaSpace{}%
\AgdaSymbol{:}\AgdaSpace{}%
\AgdaField{Objᴰ}\AgdaSpace{}%
\AgdaBound{I}\AgdaSymbol{\}}\AgdaSpace{}%
\AgdaSymbol{→}\AgdaSpace{}%
\AgdaField{Homᴰ}\AgdaSpace{}%
\AgdaSymbol{(}\AgdaFunction{C.id}\AgdaSpace{}%
\AgdaSymbol{\{}\AgdaBound{I}\AgdaSymbol{\})}\AgdaSpace{}%
\AgdaBound{x}\AgdaSpace{}%
\AgdaBound{x}\<%
\\
\>[4]\AgdaField{homomorphism}\AgdaSpace{}%
\AgdaSymbol{:}\AgdaSpace{}%
\AgdaSymbol{∀}%
\>[22]\AgdaSymbol{\{}\AgdaBound{I}\AgdaSpace{}%
\AgdaBound{J}\AgdaSpace{}%
\AgdaBound{K}\AgdaSymbol{\}}\AgdaSpace{}%
\AgdaSymbol{\{}\AgdaBound{j}\AgdaSpace{}%
\AgdaSymbol{:}\AgdaSpace{}%
\AgdaBound{C}\AgdaSpace{}%
\AgdaOperator{\AgdaFunction{[}}\AgdaSpace{}%
\AgdaBound{I}\AgdaSpace{}%
\AgdaOperator{\AgdaFunction{,}}\AgdaSpace{}%
\AgdaBound{J}\AgdaSpace{}%
\AgdaOperator{\AgdaFunction{]}}\AgdaSymbol{\}}\AgdaSpace{}%
\AgdaSymbol{\{}\AgdaBound{k}\AgdaSpace{}%
\AgdaSymbol{:}\AgdaSpace{}%
\AgdaBound{C}\AgdaSpace{}%
\AgdaOperator{\AgdaFunction{[}}\AgdaSpace{}%
\AgdaBound{J}\AgdaSpace{}%
\AgdaOperator{\AgdaFunction{,}}\AgdaSpace{}%
\AgdaBound{K}\AgdaSpace{}%
\AgdaOperator{\AgdaFunction{]}}\AgdaSymbol{\}}\AgdaSpace{}%
\AgdaSymbol{\{}\AgdaBound{x}\AgdaSymbol{\}\{}\AgdaBound{y}\AgdaSymbol{\}\{}\AgdaBound{z}\AgdaSymbol{\}}\AgdaSpace{}%
\AgdaSymbol{→}\<%
\\
\>[4][@{}l@{\AgdaIndent{0}}]%
\>[8]\AgdaField{Homᴰ}\AgdaSpace{}%
\AgdaBound{j}\AgdaSpace{}%
\AgdaBound{x}\AgdaSpace{}%
\AgdaBound{y}\AgdaSpace{}%
\AgdaSymbol{→}\AgdaSpace{}%
\AgdaField{Homᴰ}\AgdaSpace{}%
\AgdaBound{k}\AgdaSpace{}%
\AgdaBound{y}\AgdaSpace{}%
\AgdaBound{z}\AgdaSpace{}%
\AgdaSymbol{→}\AgdaSpace{}%
\AgdaField{Homᴰ}\AgdaSpace{}%
\AgdaSymbol{(}\AgdaBound{C}\AgdaSpace{}%
\AgdaOperator{\AgdaFunction{[}}\AgdaSpace{}%
\AgdaBound{k}\AgdaSpace{}%
\AgdaOperator{\AgdaFunction{∘}}\AgdaSpace{}%
\AgdaBound{j}\AgdaSpace{}%
\AgdaOperator{\AgdaFunction{]}}\AgdaSymbol{)}\AgdaSpace{}%
\AgdaBound{x}\AgdaSpace{}%
\AgdaBound{z}\<%
\\
\\[\AgdaEmptyExtraSkip]%
\>[4]\AgdaField{coe}\AgdaSpace{}%
\AgdaSymbol{:}\AgdaSpace{}%
\AgdaSymbol{∀}\AgdaSpace{}%
\AgdaSymbol{\{}\AgdaBound{I}\AgdaSpace{}%
\AgdaBound{J}\AgdaSymbol{\}}\AgdaSpace{}%
\AgdaSymbol{→}\AgdaSpace{}%
\AgdaBound{I}\AgdaSpace{}%
\AgdaOperator{\AgdaRecord{≅}}\AgdaSpace{}%
\AgdaBound{J}\AgdaSpace{}%
\AgdaSymbol{→}\AgdaSpace{}%
\AgdaField{Objᴰ}\AgdaSpace{}%
\AgdaBound{I}\AgdaSpace{}%
\AgdaSymbol{→}\AgdaSpace{}%
\AgdaField{Objᴰ}\AgdaSpace{}%
\AgdaBound{J}\<%
\\
\>[4]\AgdaField{coe-id}\AgdaSpace{}%
\AgdaSymbol{:}\AgdaSpace{}%
\AgdaSymbol{∀}\AgdaSpace{}%
\AgdaSymbol{\{}\AgdaBound{I}\AgdaSymbol{\}}\AgdaSpace{}%
\AgdaSymbol{\{}\AgdaBound{x}\AgdaSpace{}%
\AgdaSymbol{:}\AgdaSpace{}%
\AgdaField{Objᴰ}\AgdaSpace{}%
\AgdaBound{I}\AgdaSymbol{\}}\AgdaSpace{}%
\AgdaSymbol{→}\AgdaSpace{}%
\AgdaField{coe}\AgdaSpace{}%
\AgdaSymbol{(}\AgdaFunction{≅.refl}\AgdaSpace{}%
\AgdaSymbol{\{}\AgdaBound{I}\AgdaSymbol{\})}\AgdaSpace{}%
\AgdaBound{x}\AgdaSpace{}%
\AgdaOperator{\AgdaField{≈ᴰ}}\AgdaSpace{}%
\AgdaBound{x}\<%
\\
\>[4]\AgdaField{coe-comp}\AgdaSpace{}%
\AgdaSymbol{:}\AgdaSpace{}%
\AgdaSymbol{∀}\AgdaSpace{}%
\AgdaSymbol{\{}\AgdaBound{I}\AgdaSpace{}%
\AgdaBound{J}\AgdaSpace{}%
\AgdaBound{K}\AgdaSymbol{\}}\AgdaSpace{}%
\AgdaSymbol{\{}\AgdaBound{I≅J}\AgdaSpace{}%
\AgdaSymbol{:}\AgdaSpace{}%
\AgdaBound{I}\AgdaSpace{}%
\AgdaOperator{\AgdaRecord{≅}}\AgdaSpace{}%
\AgdaBound{J}\AgdaSymbol{\}}\AgdaSpace{}%
\AgdaSymbol{\{}\AgdaBound{J≅K}\AgdaSpace{}%
\AgdaSymbol{:}\AgdaSpace{}%
\AgdaBound{J}\AgdaSpace{}%
\AgdaOperator{\AgdaRecord{≅}}\AgdaSpace{}%
\AgdaBound{K}\AgdaSymbol{\}}\AgdaSpace{}%
\AgdaSymbol{\{}\AgdaBound{x}\AgdaSpace{}%
\AgdaSymbol{:}\AgdaSpace{}%
\AgdaField{Objᴰ}\AgdaSpace{}%
\AgdaBound{I}\AgdaSymbol{\}}\AgdaSpace{}%
\AgdaSymbol{→}\<%
\\
\>[4][@{}l@{\AgdaIndent{0}}]%
\>[8]\AgdaField{coe}\AgdaSpace{}%
\AgdaSymbol{(}\AgdaFunction{≅.trans}\AgdaSpace{}%
\AgdaBound{I≅J}\AgdaSpace{}%
\AgdaBound{J≅K}\AgdaSymbol{)}\AgdaSpace{}%
\AgdaBound{x}\AgdaSpace{}%
\AgdaOperator{\AgdaField{≈ᴰ}}\AgdaSpace{}%
\AgdaField{coe}\AgdaSpace{}%
\AgdaBound{J≅K}\AgdaSpace{}%
\AgdaSymbol{(}\AgdaField{coe}\AgdaSpace{}%
\AgdaBound{I≅J}\AgdaSpace{}%
\AgdaBound{x}\AgdaSymbol{)}\<%
\\
\>[4]\AgdaField{coh}%
\>[206I]\AgdaSymbol{:}\AgdaSpace{}%
\AgdaSymbol{∀}\AgdaSpace{}%
\AgdaSymbol{\{}\AgdaBound{I}\AgdaSpace{}%
\AgdaBound{J}\AgdaSymbol{\}}\AgdaSpace{}%
\AgdaSymbol{\{}\AgdaBound{x}\AgdaSpace{}%
\AgdaSymbol{:}\AgdaSpace{}%
\AgdaField{Objᴰ}\AgdaSpace{}%
\AgdaBound{I}\AgdaSymbol{\}}\AgdaSpace{}%
\AgdaSymbol{\{}\AgdaBound{I≅J}\AgdaSpace{}%
\AgdaSymbol{:}\AgdaSpace{}%
\AgdaBound{I}\AgdaSpace{}%
\AgdaOperator{\AgdaRecord{≅}}\AgdaSpace{}%
\AgdaBound{J}\AgdaSymbol{\}}\AgdaSpace{}%
\AgdaSymbol{→}\<%
\\
\>[.][@{}l@{}]\<[206I]%
\>[8]\AgdaField{Homᴰ}\AgdaSpace{}%
\AgdaSymbol{(}\AgdaField{\AgdaUnderscore{}≅\AgdaUnderscore{}.from}\AgdaSpace{}%
\AgdaBound{I≅J}\AgdaSymbol{)}\AgdaSpace{}%
\AgdaBound{x}\AgdaSpace{}%
\AgdaSymbol{(}\AgdaField{coe}\AgdaSpace{}%
\AgdaBound{I≅J}\AgdaSpace{}%
\AgdaBound{x}\AgdaSymbol{)}\<%
\\
\\[\AgdaEmptyExtraSkip]%
\>[4]\AgdaField{equiv}\AgdaSpace{}%
\AgdaSymbol{:}\AgdaSpace{}%
\AgdaSymbol{∀}\AgdaSpace{}%
\AgdaSymbol{\{}\AgdaBound{I}\AgdaSymbol{\}}\AgdaSpace{}%
\AgdaSymbol{→}\AgdaSpace{}%
\AgdaRecord{IsEquivalence}\AgdaSpace{}%
\AgdaSymbol{(}\AgdaOperator{\AgdaField{\AgdaUnderscore{}≈ᴰ\AgdaUnderscore{}}}\AgdaSpace{}%
\AgdaSymbol{\{}\AgdaBound{I}\AgdaSymbol{\})}\<%
\\
\>[4]\AgdaField{coe-resp-≈ᴰ}\AgdaSpace{}%
\AgdaSymbol{:}\AgdaSpace{}%
\AgdaSymbol{∀}\AgdaSpace{}%
\AgdaSymbol{\{}\AgdaBound{I}\AgdaSpace{}%
\AgdaBound{J}\AgdaSymbol{\}\{}\AgdaBound{x}\AgdaSpace{}%
\AgdaBound{x'}\AgdaSpace{}%
\AgdaSymbol{:}\AgdaSpace{}%
\AgdaField{Objᴰ}\AgdaSpace{}%
\AgdaBound{I}\AgdaSymbol{\}\{}\AgdaBound{I≅J}\AgdaSpace{}%
\AgdaSymbol{:}\AgdaSpace{}%
\AgdaBound{I}\AgdaSpace{}%
\AgdaOperator{\AgdaRecord{≅}}\AgdaSpace{}%
\AgdaBound{J}\AgdaSymbol{\}}\AgdaSpace{}%
\AgdaSymbol{→}\<%
\\
\>[4][@{}l@{\AgdaIndent{0}}]%
\>[8]\AgdaBound{x}\AgdaSpace{}%
\AgdaOperator{\AgdaField{≈ᴰ}}\AgdaSpace{}%
\AgdaBound{x'}\AgdaSpace{}%
\AgdaSymbol{→}\AgdaSpace{}%
\AgdaField{coe}\AgdaSpace{}%
\AgdaBound{I≅J}\AgdaSpace{}%
\AgdaBound{x}\AgdaSpace{}%
\AgdaOperator{\AgdaField{≈ᴰ}}\AgdaSpace{}%
\AgdaField{coe}\AgdaSpace{}%
\AgdaBound{I≅J}\AgdaSpace{}%
\AgdaBound{x'}\<%
\\
\>[0]\<%
\end{code}
\end{formal}

\clearpage
\subsection{Paranatural.Core}

\begin{code}[hide]%
\>[0]\AgdaSymbol{\{-\#}\AgdaSpace{}%
\AgdaKeyword{OPTIONS}\AgdaSpace{}%
\AgdaPragma{--without-K}\AgdaSpace{}%
\AgdaPragma{--safe}\AgdaSpace{}%
\AgdaSymbol{\#-\}}\<%
\\
\\[\AgdaEmptyExtraSkip]%
\>[0]\AgdaKeyword{module}\AgdaSpace{}%
\AgdaModule{Paranatural.Core}\AgdaSpace{}%
\AgdaKeyword{where}\<%
\\
\\[\AgdaEmptyExtraSkip]%
\>[0]\AgdaKeyword{open}\AgdaSpace{}%
\AgdaKeyword{import}\AgdaSpace{}%
\AgdaModule{Level}\<%
\\
\>[0]\AgdaKeyword{open}\AgdaSpace{}%
\AgdaKeyword{import}\AgdaSpace{}%
\AgdaModule{Categories.Category}\<%
\\
\>[0]\AgdaKeyword{open}\AgdaSpace{}%
\AgdaKeyword{import}\AgdaSpace{}%
\AgdaModule{Relation.Binary.PropositionalEquality}\AgdaSpace{}%
\AgdaSymbol{as}\AgdaSpace{}%
\AgdaModule{≡}\AgdaSpace{}%
\AgdaKeyword{using}\AgdaSpace{}%
\AgdaSymbol{(}\AgdaOperator{\AgdaDatatype{\AgdaUnderscore{}≡\AgdaUnderscore{}}}\AgdaSymbol{)}\<%
\\
\>[0]\AgdaKeyword{import}\AgdaSpace{}%
\AgdaModule{Categories.Morphism}\<%
\\
\>[0]\AgdaKeyword{open}\AgdaSpace{}%
\AgdaKeyword{import}\AgdaSpace{}%
\AgdaModule{GenDifunctor.Core}\<%
\\
\\[\AgdaEmptyExtraSkip]%
\>[0]\AgdaKeyword{private}\<%
\\
\>[0][@{}l@{\AgdaIndent{0}}]%
\>[2]\AgdaKeyword{variable}\<%
\\
\>[2][@{}l@{\AgdaIndent{0}}]%
\>[4]\AgdaGeneralizable{o}\AgdaSpace{}%
\AgdaGeneralizable{ℓ}\AgdaSpace{}%
\AgdaGeneralizable{e}\AgdaSpace{}%
\AgdaGeneralizable{o′}\AgdaSpace{}%
\AgdaGeneralizable{ℓ′}\AgdaSpace{}%
\AgdaGeneralizable{e′}\AgdaSpace{}%
\AgdaSymbol{:}\AgdaSpace{}%
\AgdaPostulate{Level}\<%
\\
\>[4]\AgdaGeneralizable{C}\AgdaSpace{}%
\AgdaSymbol{:}\AgdaSpace{}%
\AgdaRecord{Category}\AgdaSpace{}%
\AgdaGeneralizable{o}\AgdaSpace{}%
\AgdaGeneralizable{ℓ}\AgdaSpace{}%
\AgdaGeneralizable{e}\<%
\end{code}
\begin{formal}[\cref{defn:ParanaturalTransformation}]\label{agda:ParanaturalTransformation}\ 

\begin{code}%
\>[0]\AgdaKeyword{record}\AgdaSpace{}%
\AgdaRecord{ParanaturalTransformation}\<%
\\
\>[0][@{}l@{\AgdaIndent{0}}]%
\>[4]\AgdaSymbol{\{}\AgdaBound{C}\AgdaSpace{}%
\AgdaSymbol{:}\AgdaSpace{}%
\AgdaRecord{Category}\AgdaSpace{}%
\AgdaGeneralizable{o}\AgdaSpace{}%
\AgdaGeneralizable{ℓ}\AgdaSpace{}%
\AgdaGeneralizable{e}\AgdaSymbol{\}}\AgdaSpace{}%
\AgdaSymbol{(}\AgdaBound{D}\AgdaSpace{}%
\AgdaBound{E}\AgdaSpace{}%
\AgdaSymbol{:}\AgdaSpace{}%
\AgdaRecord{GenDifunctor}\AgdaSpace{}%
\AgdaBound{C}\AgdaSpace{}%
\AgdaGeneralizable{o′}\AgdaSpace{}%
\AgdaGeneralizable{ℓ′}\AgdaSpace{}%
\AgdaGeneralizable{e′}\AgdaSymbol{)}\<%
\\
\>[4]\AgdaSymbol{:}\AgdaSpace{}%
\AgdaPrimitive{Set}\AgdaSpace{}%
\AgdaSymbol{(}\AgdaBound{o}\AgdaSpace{}%
\AgdaOperator{\AgdaPrimitive{⊔}}\AgdaSpace{}%
\AgdaBound{ℓ}\AgdaSpace{}%
\AgdaOperator{\AgdaPrimitive{⊔}}\AgdaSpace{}%
\AgdaBound{e}\AgdaSpace{}%
\AgdaOperator{\AgdaPrimitive{⊔}}\AgdaSpace{}%
\AgdaBound{o′}\AgdaSpace{}%
\AgdaOperator{\AgdaPrimitive{⊔}}\AgdaSpace{}%
\AgdaBound{ℓ′}\AgdaSpace{}%
\AgdaOperator{\AgdaPrimitive{⊔}}\AgdaSpace{}%
\AgdaBound{e′}\AgdaSymbol{)}\AgdaSpace{}%
\AgdaKeyword{where}\<%
\\
\>[0][@{}l@{\AgdaIndent{0}}]%
\>[2]\AgdaKeyword{eta-equality}\<%
\\
\>[2]\AgdaKeyword{private}\<%
\\
\>[2][@{}l@{\AgdaIndent{0}}]%
\>[4]\AgdaKeyword{module}\AgdaSpace{}%
\AgdaModule{C}\AgdaSpace{}%
\AgdaSymbol{=}\AgdaSpace{}%
\AgdaModule{Category}\AgdaSpace{}%
\AgdaBound{C}\<%
\\
\>[4]\AgdaKeyword{module}\AgdaSpace{}%
\AgdaModule{D}\AgdaSpace{}%
\AgdaSymbol{=}\AgdaSpace{}%
\AgdaModule{GenDifunctor}\AgdaSpace{}%
\AgdaBound{D}\<%
\\
\>[4]\AgdaKeyword{module}\AgdaSpace{}%
\AgdaModule{E}\AgdaSpace{}%
\AgdaSymbol{=}\AgdaSpace{}%
\AgdaModule{GenDifunctor}\AgdaSpace{}%
\AgdaBound{E}\<%
\\
\>[4]\AgdaKeyword{open}\AgdaSpace{}%
\AgdaModule{Categories.Morphism}\AgdaSpace{}%
\AgdaBound{C}\<%
\\
\\[\AgdaEmptyExtraSkip]%
\>[2]\AgdaKeyword{field}\<%
\\
\>[2][@{}l@{\AgdaIndent{0}}]%
\>[4]\AgdaField{α}%
\>[7]\AgdaSymbol{:}\AgdaSpace{}%
\AgdaSymbol{∀}\AgdaSpace{}%
\AgdaBound{I}\AgdaSpace{}%
\AgdaSymbol{→}\AgdaSpace{}%
\AgdaFunction{D.Objᴰ}\AgdaSpace{}%
\AgdaBound{I}\AgdaSpace{}%
\AgdaSymbol{→}\AgdaSpace{}%
\AgdaFunction{E.Objᴰ}\AgdaSpace{}%
\AgdaBound{I}\<%
\\
\>[4]\AgdaField{preserve-hom}\AgdaSpace{}%
\AgdaSymbol{:}\AgdaSpace{}%
\AgdaSymbol{∀}\AgdaSpace{}%
\AgdaSymbol{\{}\AgdaBound{I}\AgdaSpace{}%
\AgdaBound{J}\AgdaSymbol{\}}\AgdaSpace{}%
\AgdaSymbol{(}\AgdaBound{j}\AgdaSpace{}%
\AgdaSymbol{:}\AgdaSpace{}%
\AgdaBound{C}\AgdaSpace{}%
\AgdaOperator{\AgdaFunction{[}}\AgdaSpace{}%
\AgdaBound{I}\AgdaSpace{}%
\AgdaOperator{\AgdaFunction{,}}\AgdaSpace{}%
\AgdaBound{J}\AgdaSpace{}%
\AgdaOperator{\AgdaFunction{]}}\AgdaSymbol{)}\AgdaSpace{}%
\AgdaSymbol{\{}\AgdaBound{x}\AgdaSymbol{\}}\AgdaSpace{}%
\AgdaSymbol{\{}\AgdaBound{y}\AgdaSymbol{\}}\AgdaSpace{}%
\AgdaSymbol{→}\<%
\\
\>[4][@{}l@{\AgdaIndent{0}}]%
\>[12]\AgdaFunction{D.Homᴰ}\AgdaSpace{}%
\AgdaBound{j}\AgdaSpace{}%
\AgdaBound{x}\AgdaSpace{}%
\AgdaBound{y}\AgdaSpace{}%
\AgdaSymbol{→}\AgdaSpace{}%
\AgdaFunction{E.Homᴰ}\AgdaSpace{}%
\AgdaBound{j}\AgdaSpace{}%
\AgdaSymbol{(}\AgdaField{α}\AgdaSpace{}%
\AgdaBound{I}\AgdaSpace{}%
\AgdaBound{x}\AgdaSymbol{)}\AgdaSpace{}%
\AgdaSymbol{(}\AgdaField{α}\AgdaSpace{}%
\AgdaBound{J}\AgdaSpace{}%
\AgdaBound{y}\AgdaSymbol{)}\<%
\\
\>[4]\AgdaField{preserve-coe}\AgdaSpace{}%
\AgdaSymbol{:}\AgdaSpace{}%
\AgdaSymbol{∀}\AgdaSpace{}%
\AgdaSymbol{\{}\AgdaBound{I}\AgdaSpace{}%
\AgdaBound{J}\AgdaSymbol{\}}\AgdaSpace{}%
\AgdaSymbol{(}\AgdaBound{I≅J}\AgdaSpace{}%
\AgdaSymbol{:}\AgdaSpace{}%
\AgdaBound{I}\AgdaSpace{}%
\AgdaOperator{\AgdaRecord{≅}}\AgdaSpace{}%
\AgdaBound{J}\AgdaSymbol{)}\AgdaSpace{}%
\AgdaSymbol{\{}\AgdaBound{x}\AgdaSpace{}%
\AgdaSymbol{:}\AgdaSpace{}%
\AgdaFunction{D.Objᴰ}\AgdaSpace{}%
\AgdaBound{I}\AgdaSymbol{\}}\AgdaSpace{}%
\AgdaSymbol{→}\<%
\\
\>[4][@{}l@{\AgdaIndent{0}}]%
\>[12]\AgdaSymbol{(}\AgdaField{α}\AgdaSpace{}%
\AgdaBound{J}\AgdaSpace{}%
\AgdaSymbol{(}\AgdaFunction{D.coe}\AgdaSpace{}%
\AgdaBound{I≅J}\AgdaSpace{}%
\AgdaBound{x}\AgdaSymbol{))}\AgdaSpace{}%
\AgdaOperator{\AgdaFunction{E.≈ᴰ}}%
\>[38]\AgdaSymbol{(}\AgdaFunction{E.coe}\AgdaSpace{}%
\AgdaBound{I≅J}\AgdaSpace{}%
\AgdaSymbol{(}\AgdaField{α}\AgdaSpace{}%
\AgdaBound{I}\AgdaSpace{}%
\AgdaBound{x}\AgdaSymbol{))}\<%
\end{code}
\end{formal}\clearpage
\subsection{StrongDinatural.Core}

\begin{code}[hide]%
\>[0]\AgdaSymbol{\{-\#}\AgdaSpace{}%
\AgdaKeyword{OPTIONS}\AgdaSpace{}%
\AgdaPragma{--without-K}\AgdaSpace{}%
\AgdaPragma{--safe}\AgdaSpace{}%
\AgdaSymbol{\#-\}}\<%
\\
\\[\AgdaEmptyExtraSkip]%
\>[0]\AgdaKeyword{module}\AgdaSpace{}%
\AgdaModule{StrongDinatural.Core}\AgdaSpace{}%
\AgdaKeyword{where}\<%
\\
\\[\AgdaEmptyExtraSkip]%
\>[0]\AgdaKeyword{open}\AgdaSpace{}%
\AgdaKeyword{import}\AgdaSpace{}%
\AgdaModule{Level}\<%
\\
\>[0]\AgdaKeyword{open}\AgdaSpace{}%
\AgdaKeyword{import}\AgdaSpace{}%
\AgdaModule{Data.Product}\<%
\\
\>[0]\AgdaKeyword{open}\AgdaSpace{}%
\AgdaKeyword{import}\AgdaSpace{}%
\AgdaModule{Categories.Category}\<%
\\
\>[0]\AgdaKeyword{open}\AgdaSpace{}%
\AgdaKeyword{import}\AgdaSpace{}%
\AgdaModule{Categories.Category.Instance.Sets}\<%
\\
\>[0]\AgdaKeyword{open}\AgdaSpace{}%
\AgdaKeyword{import}\AgdaSpace{}%
\AgdaModule{Categories.Functor}\<%
\\
\>[0]\AgdaKeyword{open}\AgdaSpace{}%
\AgdaKeyword{import}\AgdaSpace{}%
\AgdaModule{Categories.Functor.Bifunctor}\<%
\\
\>[0]\AgdaKeyword{open}\AgdaSpace{}%
\AgdaKeyword{import}\AgdaSpace{}%
\AgdaModule{Categories.Category.Product}\<%
\\
\>[0]\AgdaKeyword{open}\AgdaSpace{}%
\AgdaKeyword{import}\AgdaSpace{}%
\AgdaModule{Relation.Binary.PropositionalEquality}\AgdaSpace{}%
\AgdaSymbol{as}\AgdaSpace{}%
\AgdaModule{≡}\AgdaSpace{}%
\AgdaKeyword{using}\AgdaSpace{}%
\AgdaSymbol{(}\AgdaOperator{\AgdaDatatype{\AgdaUnderscore{}≡\AgdaUnderscore{}}}\AgdaSymbol{)}\<%
\\
\>[0]\AgdaComment{--\ open\ import\ Relation.Binary\ using\ (REL)}\<%
\\
\>[0]\AgdaComment{--\ import\ Categories.Morphism}\<%
\\
\\[\AgdaEmptyExtraSkip]%
\>[0]\AgdaKeyword{private}\<%
\\
\>[0][@{}l@{\AgdaIndent{0}}]%
\>[2]\AgdaKeyword{variable}\<%
\\
\>[2][@{}l@{\AgdaIndent{0}}]%
\>[4]\AgdaGeneralizable{o}\AgdaSpace{}%
\AgdaGeneralizable{ℓ}\AgdaSpace{}%
\AgdaGeneralizable{e}\AgdaSpace{}%
\AgdaGeneralizable{o′}\AgdaSpace{}%
\AgdaGeneralizable{ℓ′}\AgdaSpace{}%
\AgdaSymbol{:}\AgdaSpace{}%
\AgdaPostulate{Level}\<%
\\
\>[4]\AgdaGeneralizable{C}\AgdaSpace{}%
\AgdaSymbol{:}\AgdaSpace{}%
\AgdaRecord{Category}\AgdaSpace{}%
\AgdaGeneralizable{o}\AgdaSpace{}%
\AgdaGeneralizable{ℓ}\AgdaSpace{}%
\AgdaGeneralizable{e}\<%
\end{code}
\begin{formal}[\cref{defn:StrongDinaturalTransformation}]\label{agda:StrongDinaturalTransformation}\ 

\begin{code}%
\>[0]\AgdaKeyword{record}\AgdaSpace{}%
\AgdaRecord{StrongDinaturalTransformation}\AgdaSpace{}%
\AgdaSymbol{\{}\AgdaBound{C}\AgdaSpace{}%
\AgdaSymbol{:}\AgdaSpace{}%
\AgdaRecord{Category}\AgdaSpace{}%
\AgdaGeneralizable{o}\AgdaSpace{}%
\AgdaGeneralizable{ℓ}\AgdaSpace{}%
\AgdaGeneralizable{e}\AgdaSymbol{\}}\<%
\\
\>[0][@{}l@{\AgdaIndent{0}}]%
\>[4]\AgdaSymbol{(}\AgdaBound{Γ}\AgdaSpace{}%
\AgdaBound{Δ}\AgdaSpace{}%
\AgdaSymbol{:}\AgdaSpace{}%
\AgdaFunction{Bifunctor}\AgdaSpace{}%
\AgdaSymbol{(}\AgdaFunction{Category.op}\AgdaSpace{}%
\AgdaBound{C}\AgdaSymbol{)}\AgdaSpace{}%
\AgdaBound{C}\AgdaSpace{}%
\AgdaSymbol{(}\AgdaFunction{Sets}\AgdaSpace{}%
\AgdaGeneralizable{o′}\AgdaSymbol{))}\<%
\\
\>[4]\AgdaSymbol{:}\AgdaSpace{}%
\AgdaPrimitive{Set}\AgdaSpace{}%
\AgdaSymbol{(}\AgdaBound{o}\AgdaSpace{}%
\AgdaOperator{\AgdaPrimitive{⊔}}\AgdaSpace{}%
\AgdaBound{ℓ}\AgdaSpace{}%
\AgdaOperator{\AgdaPrimitive{⊔}}\AgdaSpace{}%
\AgdaBound{o′}\AgdaSymbol{)}\AgdaSpace{}%
\AgdaKeyword{where}\<%
\\
\>[0][@{}l@{\AgdaIndent{0}}]%
\>[2]\AgdaKeyword{eta-equality}\<%
\\
\>[2]\AgdaKeyword{private}\<%
\\
\>[2][@{}l@{\AgdaIndent{0}}]%
\>[4]\AgdaKeyword{module}\AgdaSpace{}%
\AgdaModule{C}\AgdaSpace{}%
\AgdaSymbol{=}\AgdaSpace{}%
\AgdaModule{Category}\AgdaSpace{}%
\AgdaBound{C}\<%
\\
\>[4]\AgdaKeyword{module}\AgdaSpace{}%
\AgdaModule{Γ}\AgdaSpace{}%
\AgdaSymbol{=}\AgdaSpace{}%
\AgdaModule{Functor}\AgdaSpace{}%
\AgdaBound{Γ}\<%
\\
\>[4]\AgdaKeyword{module}\AgdaSpace{}%
\AgdaModule{Δ}\AgdaSpace{}%
\AgdaSymbol{=}\AgdaSpace{}%
\AgdaModule{Functor}\AgdaSpace{}%
\AgdaBound{Δ}\<%
\\
\\[\AgdaEmptyExtraSkip]%
\>[2]\AgdaKeyword{field}\<%
\\
\>[2][@{}l@{\AgdaIndent{0}}]%
\>[4]\AgdaField{α}\AgdaSpace{}%
\AgdaSymbol{:}\AgdaSpace{}%
\AgdaSymbol{∀}\AgdaSpace{}%
\AgdaBound{I}\AgdaSpace{}%
\AgdaSymbol{→}%
\>[15]\AgdaFunction{Γ.F₀}\AgdaSpace{}%
\AgdaSymbol{(}\AgdaBound{I}\AgdaSpace{}%
\AgdaOperator{\AgdaInductiveConstructor{,}}\AgdaSpace{}%
\AgdaBound{I}\AgdaSymbol{)}\AgdaSpace{}%
\AgdaSymbol{→}%
\>[31]\AgdaFunction{Δ.F₀}\AgdaSpace{}%
\AgdaSymbol{(}\AgdaBound{I}\AgdaSpace{}%
\AgdaOperator{\AgdaInductiveConstructor{,}}\AgdaSpace{}%
\AgdaBound{I}\AgdaSymbol{)}\<%
\\
\>[4]\AgdaField{commute}%
\>[82I]\AgdaSymbol{:}\AgdaSpace{}%
\AgdaSymbol{∀}\AgdaSpace{}%
\AgdaSymbol{\{}\AgdaBound{I}\AgdaSpace{}%
\AgdaBound{J}\AgdaSymbol{\}}\AgdaSpace{}%
\AgdaSymbol{(}\AgdaBound{j}\AgdaSpace{}%
\AgdaSymbol{:}\AgdaSpace{}%
\AgdaBound{C}\AgdaSpace{}%
\AgdaOperator{\AgdaFunction{[}}\AgdaSpace{}%
\AgdaBound{I}\AgdaSpace{}%
\AgdaOperator{\AgdaFunction{,}}\AgdaSpace{}%
\AgdaBound{J}\AgdaSpace{}%
\AgdaOperator{\AgdaFunction{]}}\AgdaSymbol{)}\AgdaSpace{}%
\AgdaSymbol{\{}\AgdaBound{x}\AgdaSymbol{\}}\AgdaSpace{}%
\AgdaSymbol{\{}\AgdaBound{y}\AgdaSymbol{\}}\AgdaSpace{}%
\AgdaSymbol{→}\<%
\\
\>[.][@{}l@{}]\<[82I]%
\>[12]\AgdaFunction{Γ.F₁}\AgdaSpace{}%
\AgdaSymbol{(}\AgdaSpace{}%
\AgdaFunction{C.id}\AgdaSpace{}%
\AgdaOperator{\AgdaInductiveConstructor{,}}\AgdaSpace{}%
\AgdaBound{j}\AgdaSpace{}%
\AgdaSymbol{)}\AgdaSpace{}%
\AgdaBound{x}\AgdaSpace{}%
\AgdaOperator{\AgdaDatatype{≡}}\AgdaSpace{}%
\AgdaFunction{Γ.F₁}\AgdaSpace{}%
\AgdaSymbol{(}\AgdaSpace{}%
\AgdaBound{j}\AgdaSpace{}%
\AgdaOperator{\AgdaInductiveConstructor{,}}\AgdaSpace{}%
\AgdaFunction{C.id}\AgdaSymbol{)}\AgdaSpace{}%
\AgdaBound{y}\AgdaSpace{}%
\AgdaSymbol{→}\<%
\\
\>[12]\AgdaFunction{Δ.F₁}\AgdaSpace{}%
\AgdaSymbol{(}\AgdaSpace{}%
\AgdaFunction{C.id}\AgdaSpace{}%
\AgdaOperator{\AgdaInductiveConstructor{,}}\AgdaSpace{}%
\AgdaBound{j}\AgdaSpace{}%
\AgdaSymbol{)}\AgdaSpace{}%
\AgdaSymbol{(}\AgdaField{α}\AgdaSpace{}%
\AgdaBound{I}\AgdaSpace{}%
\AgdaBound{x}\AgdaSymbol{)}\AgdaSpace{}%
\AgdaOperator{\AgdaDatatype{≡}}\AgdaSpace{}%
\AgdaFunction{Δ.F₁}\AgdaSpace{}%
\AgdaSymbol{(}\AgdaSpace{}%
\AgdaBound{j}\AgdaSpace{}%
\AgdaOperator{\AgdaInductiveConstructor{,}}\AgdaSpace{}%
\AgdaFunction{C.id}\AgdaSpace{}%
\AgdaSymbol{)}\AgdaSpace{}%
\AgdaSymbol{(}\AgdaField{α}\AgdaSpace{}%
\AgdaBound{J}\AgdaSpace{}%
\AgdaBound{y}\AgdaSymbol{)}\<%
\\
\>[0]\<%
\end{code}
\end{formal}\clearpage

\end{document}